\documentclass{lms}

\usepackage{amsmath,amssymb,amscd,verbatim,graphicx}
\usepackage[all]{xy}
\SelectTips{cm}{}

\def\l{\ell}
\def\({\left(}
\def\too{\longrightarrow}
\def\){\right)}
\def\[{\psi\(\tfrac12\<}
\def\]{\>\)}

\def\inv{^{-1}}

\def\dens#1#2{\left|\det #1\right|^{#2}}
\def\lll{_{\l_1,\l_2}}

\def\gll{_{g\l,\l}}

\DeclareMathOperator\supp{support}

\def\Tr{\mathrm{Tr\,}}
\DeclareMathOperator\End{End}

\def\SSS{\mathfrak S}
\def\ooo{\mathfrak o}
\def\HHH{\mathcal H}

\def\ZZZ{\mathbb Z}
\def\FF{\mathcal F}
\def\ZZZ{\mathbb Z}

\def\Sp{{\mathit{Sp}}}

\def\Mp{{\mathit{Mp}}}
\def\SL{{\mathit{SL}}}

\def\CCC{\mathbb C}
\def\RRR{\mathbb R}

\def\GL{{\mathit{GL}}}
\def\gl{{\mathfrak{gl}}}

\def\Grass{\mathrm{Gr}}

\def\<{\left<}
\def\>{\right>}

\newtheorem{prop}{Proposition}[section]
\newtheorem{lemmalrestrict'}{Lemma \ref{lrestrict}{\rm(iii)}$'$}
\newtheorem{proprel'}{Proposition \ref{rel}$'$}

\newtheorem{lemma}[prop]{Lemma}
\newtheorem{corollary}[prop]{Corollary}

\newunnumbered{definition}{Definition}

\newunnumbered{rem}{Remark}
\newnumbered{remnum}{Remark}[section]

\def\gl{_{g,\l}}
\def\gll{_{g\l,\l}}

\title{The Character of the Weil Representation}
\date{1 January 2007}
\author{Teruji Thomas}
\classno{11F27}
\extraline{This work was partially supported by NSF grants DMS-9977134 and DMS-0401164 and by the University of Chicago.}

\begin{document}\maketitle
\begin{abstract} Let $V$ be a symplectic vector space over a finite or local field. We compute the character of the Weil representation of the metaplectic group $\Mp(V)$. The final formulas are overtly free of choices (e.g. they do not involve the usual choice of a Lagrangian subspace of $V$).  Along the way, in results similar to those of K. Maktouf, we relate the character to the Weil index of a certain quadratic form, which may be understood as a Maslov index.  This relation also expresses the character as the pullback of a certain simple function from $\Mp(V\oplus V)$. 
\end{abstract}

\section{Introduction}

 Let $F$ be a field of characteristic not $2$: it may be the real numbers $\RRR$, the complex numbers $\CCC$, a non-archimedean local field, or a finite field. Fix a non-trivial additive character $\psi\colon F\to\CCC^\times$. Let $V$ be a vector space over $F$ with symplectic form $\<-\,,-\>$.

The symplectic group $\Sp(V)$ has a well known projective representation $\rho$, de\-pending on $\psi$, called the `Weil representation' after A. Weil's seminal paper \cite{We}. When $F$ is finite or complex, $\rho$ can be realized as a true representation of $\Sp(V)$, and in the other cases as a true representation of a two-fold cover $\Mp(V)$ of $\Sp(V)$. In this article we compute the characters $\Tr\rho$ of these representations.

The standard constructions of $\rho$ involve a choice: for example, that of a Lagran\-gian subspace $\l\subset V$. Our goal, more precisely, is to present formulas for $\Tr\rho$ completely free of such choices.

\subsection{Formulation of the Main Results.}

\subsubsection{}\label{first}
 First suppose $F$ is finite of cardinality $|F|$, so $\Tr\rho$ is a function on $\Sp(V)$.

Given $g\in\Sp(V)$, the endomorphism $(g-1)\in\End(V)$ plays a key role in the formula for $\Tr\rho(g)$ (cf. \cite{GH} \S2.1 and \cite{Howe} p. 294). Let us denote by $\sigma_g$ the induced isomorphism
\begin{equation*}\sigma_g\colon V/\ker(g-1)\overset\sim\too (g-1)V.\end{equation*}
It is easy to check that $v\otimes w\mapsto\<\sigma_g v,w\>$ defines a nondegenerate bilinear form on $V/\ker(g-1)$. 
Let $\det\sigma_g\in F^\times/(F^\times)^2$ be its discriminant (see \S\ref{orient}); if $(g-1)$ is invertible then $\det\sigma_g$ is just the usual determinant $\det(g-1)$ mod $(F^\times)^2$.

The second ingredient we need is the `Weil index,' a character $\gamma$ of the Witt group $W(F)$ of quadratic forms over $F$ (see \S\ref{WittWeil} and \cite{We},  \cite{Pe}).  If $a\in F^\times$ then denote by $\gamma(a)$ the value of $\gamma$ on the one-dimensional quadratic form $x\mapsto ax^2$.  In this finite field case,
\begin{equation*}\gamma(a)=|F|^{-1/2}\sum_{x\in F}\psi(\tfrac12ax^2).\end{equation*}
It depends only on $a\bmod(F^\times)^2.$
\begin{theorem1a*}\label{THM1A} If $F$ is a finite field,  then
\begin{equation*}\Tr\rho(g)=|F|^{\tfrac12\dim\ker(g-1)}\gamma(1)^{\dim V-\dim\ker(g-1)-1}\gamma(\det\sigma_g).\end{equation*}
\end{theorem1a*}

\begin{remnum} Restricted to the case where $\ker(g-1)=0$, an equivalent formula has been independently obtained by S. Gurevich and R. Hadani \cite{GH}, as a corollary to their algebraic-geometric approach to the Weil representation over a finite field. 
\end{remnum}

\begin{remnum} \label{remHowe} R. Howe \cite{Howe} understood many aspects of Theorem 1A without, apparently, finding a closed formula. For example, one can determine the absolute value of $\Tr\rho(g)$ from the fact that $\rho\otimes\rho^*$ is the natural action of $\Sp(V)$ on $L^2(V)$.
\end{remnum}

\begin{remnum}
\label{remchiit}
The values of $\gamma$ depend on $\psi$.  For $a\in  F^\times$, $\gamma(a)$ is calculated for standard choices of $\psi$ in the appendix of \cite{Pe}.
It is well known  that, in this finite field case, $\chi\colon a\mapsto \gamma(-1)\gamma(a)$ is the unique non-trivial character of $F^\times/(F^\times)^2\cong\ZZZ/2\ZZZ$ (see e.g. \cite{Bu}, Exercise 4.1.14), and consequently $\gamma(a)\gamma(b)=\gamma(1)\gamma(ab)$.  With this in mind one can  re-write Theorem 1A in various ways, for example as\begin{equation*}\Tr\rho(g)=|F|^{\tfrac12\dim\ker(g-1)}\gamma(1)^{\dim V-\dim\ker(g-1)}\chi(\det\sigma_g).\end{equation*}
In \S\ref{example} we consider in detail what happens when $\dim V=2$.
\end{remnum}

\subsubsection{} Now suppose $F=\CCC$ is the field of complex numbers. Then $\Tr\rho$ is defined as a generalized function on $\Sp(V)$ (see \S\ref{char}). It is known by the work of Harish Chandra to be (represented by) a locally integrable function, smooth (i.e. $C^\infty$) on the open set of regular semisimple elements. 
Let 
\begin{equation}\label{Sp''}\Sp(V)'':=\{g\in\Sp(V)\,|\,\det(g-1)\neq0\}.
\end{equation} In \cite{Ma}, p. 293, it is shown that $\Sp(V)''$ contains the regular semisimple elements, but it obviously contains much more. 
\begin{theorem1b*}\label{THM1B} 
 If $F=\CCC$ then $\Tr\rho$ is smooth on $\Sp(V)''$, and indeed
\begin{equation*}\Tr\rho(g)=\dens{(g-1)}{-1}.\end{equation*}
\end{theorem1b*}

This statement appears as Theorem 1 in part II of \cite{Torasso}.
Perhaps it is worth remarking that $\gamma\equiv 1$ when $F=\CCC$. 

\subsubsection{} Finally, suppose $F$ is the field of real numbers or else a non-archimedean local field.  Now $\rho$ is a representation not of $\Sp(V)$ but of a double cover $\Mp(V)$. A standard construction of $\Mp(V)$ is recalled in \S\ref{Meta}. Let $\pi\colon \Mp(V)\to\Sp(V)$ be the projection and \begin{equation}\label{DMp''}\Mp(V)'':=\pi\inv(\Sp(V)'')=\{\tilde g\in\Mp(V)\,|\,\det(\pi(\tilde g)-1)\neq0\}.
\end{equation}

The character $\Tr\rho$ is a generalized function on $\Mp(V)$, again known to be represented by a locally integrable function, smooth  on the set of regular semisimple elements. Here $\tilde g\in\Mp(V)$ is said to be regular semisimple if $\pi(\tilde g)$ itself is;  `smooth' means `$C^\infty$' in the real and `locally constant' in the non-archimedean cases. 

If $\tilde 1\in\Mp(V)$ is the non-identity element over $1\in\Sp(V)$, then $\rho(\tilde 1)=-1$. Thus,  given $\tilde g\in\Mp(V)$, $\Tr\rho(\tilde g)$ is determined up to sign by $g:=\pi(\tilde g)\in \Sp(V)$.

\begin{theorem1c*}\label{THM1C} 
 If $F$ is real or non-archimedean, then $\Tr\rho$ is smooth on $\Mp(V)''$, and given up to sign by 
\begin{equation*}\Tr\rho(\tilde g)=\pm\frac{\gamma(1)^{\dim V-1}\gamma(\det(g-1))}{\dens{\(g-1\)}{1/2}}.\end{equation*}
\end{theorem1c*}

\begin{remnum} Various aspects of $\Tr\rho$  were previously understood, including some formulas (see e.g. \cite{Ad}, \cite{Howe}, \cite{Ma}, \cite{Torasso},    as well as Remark \ref{remHowe}  and \S\ref{mak1}). What seems to be new in Theorem 1C is that it is overtly independent of choices, as explained at the beginning of this article. The right-hand side is also easy to compute using the values of $\gamma$ from \cite{Pe}.
\end{remnum}

\subsection{The Character (Without Sign Ambiguity) via the Maslov Index.} 

In Theorem 2, stated in this section and proved beginning in
\S\ref{overview}, we express $\Tr\rho$ in terms of the Maslov index $\tau$  (see \S\ref{RMI} and \cite{LV},\cite{Th}). Alternatively, as we explain in \S\ref{zowee}, Theorem 2 can be understood to describe $\Tr\rho$ as the pullback of a simple function from a larger metaplectic group.
Either way, this computes $\Tr\rho$ with no sign ambiguity, in contrast to Theorem 1C; but to get a number, one must choose a Lagrangian subspace of $V$.

We deduce Theorem 1 from Theorem 2 in \S\ref{implication}.

\subsubsection{}\label{evaluate} To give a uniform approach, define $\Mp(V)$, when $F$ is complex or finite, as the trivial extension of $\Sp(V)$ by $\ZZZ/2\ZZZ$.  In general, we we will use the following description of $\Mp(V)$ taken from \cite{LV},\cite{Ma},\cite{Pe}. 

Let $\Grass(V)$ be the set of Lagrangian subspaces in $V$. 
An element of $\Mp(V)$ is a pair $(g,t)$ with $g\in\Sp(V)$ and $t$ a function $t\colon \Grass(V)\to \CCC^\times$ satisfying certain conditions
 (we will recall more details in \S\ref{Meta}). 

\begin{remnum}\label{splitting}
This description makes sense over any of our fields. When $F=\CCC$ the splitting of $\Mp(V)\to\Sp(V)$ is given by $g\mapsto(g,1)$. However, when $F$ is finite, the splitting is more complicated (see Proposition \ref{split}).
\end{remnum}

\subsubsection{}\label{ee}  Fix $\l\in\Grass(V)$. Let $\overline V$ be $V$ equipped with minus the given symplectic form. Let $\Gamma_g$ be the graph of $g\colon\overline V\to V$, so $\Gamma_1$ is the diagonal in $\overline V\oplus V$. Then  $\Gamma_g$, $\Gamma_1$, and $\l\oplus \l$ are all Lagrangian subspaces of $\overline V\oplus V$.  For $(g,t)\in\Mp(V)$ set
\begin{equation}\label{eqTheta2}
\Theta_\l(g,t):=t(\l)\cdot\gamma(\tau(\Gamma_g,\Gamma_1,\l\oplus \l))
\end{equation}
where $\tau$ is the Maslov index.
In \S\ref{constancy} we prove

\begin{prop}\label{indep} $\Theta_\l$ is independent of $\l$ and locally constant on $\Mp(V)''$.\end{prop}
$\Mp(V)''$ was defined by formula \eqref{DMp''}.
 \begin{remnum}   Neither factor in the definition \eqref{eqTheta2} of $\Theta_\l$ is itself continuous on $\Mp(V)''$, and both depend on the choice of $\l$. A more canonical description of $\Theta_\l$ is given in \S\ref{zowee}. \end{remnum}

\stepcounter{theorem}
\begin{theorem2a*}\label{thm2A}   Suppose $F$ is infinite. Then 
 \begin{equation*}\Tr\rho(g,t)=\left\|\det(g-1)\right\|^{-1/2}\cdot \Theta_\l(g,t)\end{equation*}
where $\|\cdot\|$ denotes the usual norm when $F$ is real or archimedean, and  the square of the usual norm when $F=\CCC$.
\end{theorem2a*}

\begin{remnum} If $F=\CCC$ then $\gamma\equiv1$ so $\Theta_\l(g,1)=1$. Thus Theorem 1B already follows from Theorem 2A and Remark \ref{splitting}.
\end{remnum}

\begin{theorem2b*}\label{thm2b} Suppose $F$ is finite. Then 
\begin{equation*}\Tr\rho(g,t)=|F|^{\tfrac12\dim\ker(g-1)}\cdot \Theta_\l(g,t).\end{equation*}
\end{theorem2b*}

\subsubsection{An Explicit Quadratic Form.}
\label{defq1}

  In \S\ref{theform} we describe a canonical quadratic space $(S'\gl,q'\gl)$ representing the Maslov index $\tau(\Gamma_g,\Gamma_1,\l\oplus\l)$ that appears in \eqref{eqTheta2}. For $g\in\Sp(V)''$ the answer is particularly simple (and this is the only case needed for Theorem 2A): 

For any fixed $\l\in\Grass(V)$ and any $g\in\Sp(V)''$ define 
\begin{equation*}\label{eqq1}q'\gl(a,b)=\<a,(g-1)\inv b\>
\end{equation*} for all $a,b\in \l$. Then $q'\gl$ is a  symmetric bilinear form on $\l$, and $S'\gl$ is defined to be the quotient of $\l$ on which 
$q'\gl$ is nondegenerate. (The symmetry of $q'\gl$ is explained in Proposition \ref{aresym}.)

\subsubsection{Relation to the Work of Maktouf.}
 \label{mak1}  K. Maktouf \cite{Ma} proved a theorem in the $p$-adic case very similar to Theorem 2A.  As explained in \S\ref{maktouf}, our formulas are identical when $g$ is semisimple and $\l$ is appropriately chosen---in general, Maktouf's version of $\Theta_\l$ (denoted $\Phi$ in \cite{Ma}), and even $\l$ itself, is constructed from the semisimple part of $g$. Since the semisimple elements are dense in $\Mp(V)$, one can deduce the $p$-adic Theorem 2A from the main theorem in \cite{Ma} and Proposition \ref{indep}. However, our proof is different and in some sense more direct.

\subsection{The Character as the Pullback of a Function from $\Mp(\overline V\oplus V)$.}\label{zowee}

We use the notation of \S\ref{evaluate}--\S\ref{ee}.
Let $f$ be the embedding $f\colon\Sp(V)\to\Sp(\overline V\oplus V)$ defined by $f(g)=(1,g)$. 
There is a unique\footnote{Unless $|F|=3$ and $\dim V=2$, when $\Mp(V)$ has automorphisms over $\Sp(V)$---see Remark \ref{splitrem}.} homomorphic embedding $\tilde f\colon\Mp(V)\to\Mp(\overline V\oplus V)$ covering $f$, described explicitly in \S\ref{constancy}. Consider the function $\mathrm{ev}_{\Gamma_1}$ on $\Mp(\overline V\oplus V)$ defined by $(g',t')\mapsto t'(\Gamma_1)$.  

\begin{prop}\label{zow} $\Theta_\l(g,t)=\mathrm{ev}_{\Gamma_1}\circ \tilde f$. 
\end{prop}

Along with Theorem 2, this gives another description of $\Tr\rho$.
The proof of Proposition \ref{zow} is given in \S\ref{constancy}, where we use it to deduce Proposition \ref{indep}. 

\begin{acknowledgements} I would like to thank J. Adams for inspiring this work, and him and R. Kottwitz for some useful comments. I am also grateful to the referee for a careful reading.\end{acknowledgements}

\section{Example: $\SL_2$ over a Finite Field}\label{example}

 Suppose $F$ is a finite field (but the infinite case is similar); let $V=F^2$ with basis $\{e_1,e_2\}$ such that $\<e_1,e_2\>=1$. Then $\Sp(V)=\SL_2(F)$. Suppose $g=\(\begin{smallmatrix} a&b\\c&d\end{smallmatrix}\),$ with $ad-bc=1$. Then $\det(g-1)=2-a-d$.  With similar calculations, and using Theorem 1A and Remark \ref{remchiit}, one finds:
\begin{enumerate}
\item[(i)] If $a+d\neq2$ then $\Tr\rho(g)=\gamma(1)^2\chi(2-a-d)=\chi(a+d-2)$.
\item[(ii)] If $a+d=2$, $b\neq 0$, then $\det\sigma_g=b\bmod(F^\times)^2$, $\Tr\rho(g)=|F|^{1/2}\gamma(1)\chi(b)$.
\item[(iii)] If $a+d=2$, $c\neq 0$, then $\det\sigma_g=-c\bmod(F^\times)^2$, $\Tr\rho(g)=|F|^{1/2}\gamma(1)\chi(-c)$.
\item[(iv)] If $a+d=2,$ $b=c=0$, then $g=1$, $\det\sigma_g=1\bmod(F^\times)^2$, $\Tr\rho(g)=|F|$.
\end{enumerate}
Let us test these formulas on some standard elements of $\Sp(V)$. 
The results may be verified using explicit formulas for $\rho$, like those in \cite{Bu}, Exercises 4.1.13--4.1.14.
\begin{enumerate}
\item[(a)] $\Tr\rho\(\begin{smallmatrix} a&b\\0&1/a\end{smallmatrix}\)=\chi(a)$ if $a\neq 1$ (note that $a+a\inv-2\equiv a\bmod(F^\times)^2$).
\item[(b)] $\Tr\rho\(\begin{smallmatrix} 1&b\\0&1\end{smallmatrix}\)=|F|^{1/2}\gamma(1)\chi(b)$ if $b\neq 0$. 
\item[(c)] $\Tr\rho\(\begin{smallmatrix} 0&1\\-1&0\end{smallmatrix}\)=\chi(-2)$. In case $|F|=p$ prime, Quadratic Reciprocity implies that
 $\chi(-2)=1$ if $p\equiv1,3\bmod 8$, and $\chi(-2)=-1$ if $p\equiv5,7\bmod 8$.
\end{enumerate}

\section{Recollections I: Weil Index}

 We call a {\em quadratic space} a vector space equipped with a nondegenerate symmetric bilinear form.  Let $W(F)$ be the Witt group (or ring, but we are only interested in the additive structure) of quadratic spaces over $F$ (see \cite{Lam}).

\label{WittWeil}

\subsection{Definition and Properties.}\label{sg}

 The next proposition/definition is due to Weil \cite{We} (see also \cite{LV} p. 111). 
 If $(A,q)$ is a quadratic space, we write $f_{q}$ for the function $x\mapsto\psi(\tfrac12q(x,x))$ on $A$; $dq$ 
is the self-dual measure on $A$,
 so that the Fourier transform $(f_{q}\,dq)^\wedge$ is a generalized function on $A^*$; and $q^*$ is the dual quadratic form on $A^*$. That is, $q^*(x,y):=\<x,\Phi^{-1}(y)\>$ where 
$\Phi\colon A\to A^*$ is the isomorphism $a\mapsto q(a,-)$. 

\begin{prop}\label{pd}There exists a number $\gamma(q)$ such that
\begin{equation}\label{eqg}
(f_{q}\,dq)^\wedge=\gamma(q)\cdot f_{-q^*}\end{equation}
  as generalized functions on $A^*$. Moreover, $(A,q)\mapsto\gamma(q)$ is a (unitary) character $W(F)\to\CCC^\times$. 
If $a\in F^\times$ then denote by $\gamma(a)$ the value of $\gamma$ on the one-dimensional quadratic form $x\mapsto ax^2$. Then:
\begin{enumerate}
\item[(i)] For $a\in F^\times$, $\gamma(a)$ depends only on $a\bmod(F^\times)^2$.
\item[(ii)] $\gamma(a)\gamma(b)=\pm\gamma(1) \gamma(ab)$ for all $a,b\in F^\times$, with a plus if $F$ is finite or complex.
\item[(iii)] For a quadratic space $(A,q)$, $\gamma(q)=\pm\gamma(1)^{\dim A-1}\gamma(\det q)$, with a plus when $F$ is finite or complex; here $\det q$ is the discriminant of $q$.
\end{enumerate}
\end{prop}

To be precise, in (ii) the sign is the Hilbert symbol $(a,b)$: $(a,b)=1$ if $a$ is a norm from $F[\sqrt b]$, and $(a,b)=-1$ otherwise.  In (iii) the sign is the Hasse invariant of $q$.

\subsection{Integral Formulas.}

First suppose $F$ is finite. We allow here the case when $q$ may be degenerate on $A$. Then to define $\gamma(q)$ one should consider $q$ as a nondegenerate form on $A/\ker q$.
\begin{prop} \label{pg1}  Let $q$ be a possibly degenerate quadratic form on $A$. Then 
\begin{equation*}\gamma(q)=|F|^{-\tfrac12(\dim A/\ker q)-\dim\ker q}\sum_{x\in A} \psi(\tfrac12q(x,x)).\end{equation*}
\end{prop}
\begin{proof} If $\ker q=0$ then the right side is just $(f_q\,dq)^\wedge$, evaluated at $0$. In general, use that $f_q$ is constant along the cosets of $\ker q$ in $A$.
\end{proof}

\subsubsection{}\label{Dirac} Now we treat $F$ infinite. Let $h$ be a Schwartz function on $A$ such that $h(0)=1$ and such that the Fourier transform $h^\wedge$ is a positive measure.   For $s\in F$ set $h_s=h(sx)$; then the Fourier transform $h_s^\wedge$ approaches a delta-measure as $s\to 0$. 

\begin{prop} \label{pg}  Suppose $F$ is infinite and $q$ nondegenerate. For $h$ as above, 
\begin{equation*}\gamma(q)=\lim_{s\to 0}\int_{x\in A} h_s(x)\cdot \psi(\tfrac12q(x,x))\,dq.\end{equation*}
\end{prop}

\begin{proof} Certainly 
$\displaystyle\int_{x\in A} h_s(x)\cdot\psi(\tfrac12q(x,x))\,dq=
\displaystyle\int_{x\in A^*} h_s^\wedge(x)\cdot(f_q\,dq)^\wedge(x).$
Now apply \eqref{eqg} and take the limit $s\to0$.
\end{proof}

\section{Recollections II: Maslov Index}  \label{RMI}

Let $U$ be a symplectic vector space (in our applications, $U=V$ or $U=\overline V\oplus V$).

\subsection{Basic Properties.}\label{mi} For any $n>0$, the Maslov index is a function \begin{equation*}\tau\colon\Grass(U)^n\to W(F),\end{equation*} invariant under the diagonal action of $\Sp(U)$ on $\Grass(U)^n$. 
We recall the following properties (for which see \cite{Th}, \cite{KS}):
\begin{enumerate}
\item[(i)] Dihedral symmetry:
\begin{equation*}\tau(\l_1,\ldots,\l_n)=-\tau(\l_n,\ldots,\l_1)=\tau(\l_n,\l_1,\ldots,\l_{n-1}).\end{equation*}
\item[(ii)] Chain condition: For any $j$, $1<j<n$,
\begin{equation*}\tau(\l_1,\l_2,\ldots, \l_j)+\tau(\l_1,\l_j,\ldots,\l_n)=\tau(\l_1,\l_2,\ldots,\l_n).\end{equation*}
\item[(iii)] Additivity: If $U,U'$ are symplectic spaces, $\l_1,\ldots,\l_n\in\Grass(U),$ $\l'_1,\ldots,\l'_n\in\Grass(U')$, 
so that $\l_i\oplus \l_i'\in\Grass(U\oplus U')$, we have
\begin{equation*}\tau(\l_1\oplus \l_1',\ldots, \l_n\oplus \l_n')=\tau(\l_1,\ldots,\l_n)+\tau(\l_1',\ldots,\l_n').\end{equation*}
\end{enumerate}

\subsection{Rank and Discriminant.} In \cite{Th} we constructed a canonical quadratic space representing $\tau(\l_1,\ldots,\l_n)$. We now will give formulas for the rank (stated in \cite{Th}) and the discriminant (essentially found in \cite{Pe}). For the latter we need the notion of an oriented Lagrangian.

\subsubsection{ Orientations and Discriminants.}

\label{orient} An oriented vector space is a pair $(A,o)$ where $A$ is a vector space and $o$ is an element of $\det A$, consided up to multiplication by $(F^\times)^2$. Write $\ooo(A)$ for the set of all orientations on $A$. 

\begin{enumerate}
\item[(i)] For any subspace $B\subset A$, the wedge product gives a natural map 
\begin{equation*}\wedge\colon \ooo(B)\times_{F^\times}\ooo(A/B)\to\ooo(A).\end{equation*} 
\item[(ii)]
 For any perfect pairing $\beta\colon A\otimes A'\to F$, there is also a pairing 
\begin{equation*}\beta\colon\ooo(A)\times_{F^\times}\ooo(A')\to F^\times/(F^\times)^2.\end{equation*}
\end{enumerate}
Namely, if $\{e_1,\ldots,e_m\}$ is a basis for $A$ and $\{e_1^*,\ldots,e_m^*\}$ is the dual basis for $A'$, then $\beta(e_1\wedge\cdots\wedge e_m,e^*_1\wedge\cdots\wedge e_m^*)=1$.

\begin{definition} In case $A=A'$, the number $\beta(o,o)\in F^\times/(F^\times)^2$ is independent of the orientation $o$ on $A$, and is called the {\em discriminant} of $\beta$. (If $A=0$, then the discriminant is defined to be $1$.)
\end{definition}

\subsubsection{}\label{DA} Suppose $(\l_1,o_1)$ and $(\l_2,o_2)$ are Lagrangians with orientations. Choose any orientation $o$ of $\l_1\cap\l_2$.  For $i=1,2$ there is a unique orientation $\bar o_i$ of $\l_i/\l_1\cap\l_2$ such that $o\wedge\bar o_i=o_i$. The symplectic form induces a perfect pairing $(\l_1/\l_1\cap\l_2)\otimes (\l_2/\l_1\cap\l_2)\to F$. 
Set \begin{equation*}O\lll:=\<\bar o_1,\bar o_2\>.\end{equation*} It is independent of $o$ but obviously depends on $o_1,o_2$.

\subsubsection{} From \cite{Th} in conjunction with the calculations in \cite{Pe} \S1.6.1, we obtain 

\begin{prop}\label{rankdisc} The Maslov index $\tau(\l_1,\ldots,\l_n)$ can be represented by a canonically defined quadratic space $(T,q)$ with dimension
\begin{equation*}\dim T=\frac{n-2}{2}\dim V-\sum_{i\in\ZZZ/n\ZZZ}\dim (\l_i\cap \l_{i+1})+2\dim \bigcap_{i\in\ZZZ/n\ZZZ} \l_i\end{equation*} and discriminant
\begin{equation*}\det {q}=
(-1)^{\tfrac12\dim V+\dim\bigcap_{i\in\ZZZ/n\ZZZ} \l_i}\prod_{i\in\ZZZ/n\ZZZ} O_{\l_i,\l_{i+1}}\end{equation*} 
for arbitrarily chosen orientations on $\l_1,\ldots,\l_n$.\end{prop}

\subsection{Weil Index.} Fix an arbitrary orientation on each $\l_i$. Define 
\begin{equation}\label{eqm}m(\l_i,\l_{i+1}):=
\gamma(1)^{\tfrac12\dim V-\dim\l_i\cap \l_{i+1}-1}\gamma(O_{\l_{i},\l_{i+1}}).\end{equation}
Applying Propositions \ref{rankdisc} and \ref{pd}, one deduces as in \cite{Pe}
\begin{corollary}\label{cor:m} 
$\gamma(\tau(\l_1,\ldots,\l_n))=\pm \prod_{i\in\ZZZ/n\ZZZ}m(\l_i,\l_{i+1})$, with a plus when $F$ is finite or complex.
\end{corollary}

\begin{remnum} On the right side of Corollary \ref{cor:m}, the factors $m(\l_i,\l_{i+1})$ depend on the orientations on $\l_i,\l_{i+1}$, but the overall product does not.
\end{remnum}

\begin{remnum}\label{remA}  Regarding the compatibility of Corollary \ref{cor:m} with \S\ref{mi}(i,ii,iii) and the symplectic invariance of $\tau$, one can easily check 
\begin{equation}\label{minv}\begin{aligned}
 &m(g\l_i,g\l_{i+1})=m(\l_i,\l_{i+1})=\pm m(\l_{i+1},\l_i)\inv\\
 &m(\l_i\oplus\l'_i,\l_{i+1}\oplus\l'_{i+1})=\pm m(\l_i,\l_{i+1})\cdot m(\l'_i,\l'_{i+1})\end{aligned}
 \end{equation}
 with plusses as usual if $F$ is finite or complex. 
\end{remnum}

\section{Recollections III: Metaplectic Group}\label{Meta}

In this section we recall the explicit model of the metaplectic group developed in \cite{Pe},\cite{LV}. Again let $U$ be a symplectic vector space; in practice $U=V$ or $U=\overline V\oplus V$.

\subsection{The Maslov Cocycle.}
 Fix $\l\in\Grass(U)$. The properties of the Maslov index recalled in \S\ref{mi}(i,ii) imply that the function
\begin{equation}\label{cgh}
(g,h)\mapsto c_{g,h}(\l):=\gamma(\tau(\l,g\l,gh\l))
\end{equation}
is a $2$-cocycle on $\Sp(U)$ with values in  $\CCC^\times$. 
Corollary \ref{cor:m} implies that one can change this cocycle by a coboundary to take values in $\{\pm 1\}\subset\CCC^\times$. The metaplectic group $\Mp(U)$  is defined to be the corresponding central extension of $\Sp(U)$ by $\{\pm1\}$. 
Let us now make this construction explicit.

\subsection{Definition of the Metaplectic Group.} \label{constr}
 Observe that, in the notation of \eqref{eqm},
\begin{equation}\label{eqs}
m_g(\l):=m(g\l,\l)=\gamma(1)^{\tfrac12\dim U-\dim g\l\cap\l-1}\gamma(O\gll)
\end{equation}
is independent of the choice of orientation on $\l$, if we give $g\l$ the same orientation transported by $g$. 

\begin{definition}
 The metaplectic group $\Mp(U)$ consists of pairs of the form $(g,\pm t_g)$ 
 with $g\in\Sp(U)$ and $t_g\colon\Grass(U)\to\CCC$ any function satisfying the following conditions:
\begin{enumerate}
\item[(i)] $t_g(\l)^2=m_g(\l)^2$ and  
\item[(ii)] $t_g(\l')=\gamma(\tau(\l,g\l,g\l',\l'))\cdot t_g(\l)$ for any $\l,\l'\in\Grass(U)$.
\end{enumerate}
Multiplication in $\Mp(U)$ is given by 
\begin{equation}\label{multiply}(g,s)\cdot(h,t)=(gh,s t c_{g,h})\end{equation}
where $c_{g,h}$ is defined by \eqref{cgh}. 
\end{definition}

This particular form of the construction appears in \cite{LP} and apparently goes back to M. Duflo.

\subsection{Properties.}

\begin{prop}\label{exactlytwo} $\Mp(U)$ is a two-fold cover of $\Sp(U)$, i.e. there are exactly two functions $t_g\colon\Grass(U)\to\CCC$ satisfying conditions (i) and (ii).
 In case $F$ is finite or complex, we may take $t_g=\pm m_g$.
\end{prop}
\begin{proof} There are clearly zero or two such functions. To show they do exist, one has to verify the following two facts for all $\l,\l',\l''\in\Grass(U)$.
\begin{enumerate}
\item[(a)] $m_g(\l')=\pm \gamma(\tau(\l,g\l,g\l',\l'))\cdot m_g(\l)$ with a plus if $F$ is finite or complex.
\item[(b)] $\tau(\l,g\l,g\l',\l')+\tau(\l',g\l',g\l'',\l'')=\tau(\l,g\l,g\l'',\l'').$
\end{enumerate}
 Use Corollary \ref{cor:m} and property \eqref{minv} to prove (a). As for (b), \S\ref{mi}(i,ii) imply
\begin{equation*} 
 \tau(\l,g\l,g\l',\l')+\tau(\l',g\l',g\l'',\l'')=\tau(\l,g\l,g\l'',\l'')+\tau(g\l,g\l',g\l'')+\tau(\l,\l'',\l').\end{equation*}
 But the last two terms cancel by \S\ref{mi}(i) and the symplectic invariance of $\tau$.
\end{proof}

\begin{prop}\label{split} If $F$ is finite or complex then $g\mapsto(g,m_g)$ is a group homomorphism splitting the projection $\Mp(U)\to\Sp(U)$. If $F$ is complex then $m_g\equiv 1$.
\end{prop}

\begin{proof} For the first statement we must check that $m_gm_hc_{gh}=m_{gh}$. This is immediate from \eqref{cgh}, \eqref{minv}, and Corollary \ref{cor:m}. The second statement is obvious since $\gamma\equiv 1$ when
 $F=\CCC$.
\end{proof}

\begin{remnum}\label{splitrem}  The splitting of Proposition \ref{split} is known to be unique except when $|F|=3$ and $\dim U=2$.
\end{remnum}

\begin{prop}\label{topology}  The topology on $\Mp(U)$ is determined by the fact that $\pi\colon\Mp(U)\to\Sp(U)$ is a local homeomorphism, and by the following property. For each $\l\in\Grass(U)$ and $n\geq 0$, set $N_{\l,n}:=\{g\in\Sp(U)\,|\,\dim \l\cap g\l=n\}$. Then the function $\mathrm{ev}_\l\colon(g,t)\mapsto t(\l)$ is locally constant over each $N_{\l,n}$.
\end{prop}
\begin{proof} That the property holds is implicit in \cite{LV}, \S1.9.11. To determine the topology on the group, it suffices to specify a neighbourhood $N$ of a single point $(g,t)$ such that $\pi$ is injective on $N$. We can choose $g$ to lie in the open set $N_{\l,0}$,  and take $N=\{(h,s)\,|\, h\in N_{\l,0},\, \mathrm{ev}_\l(h,s)=\mathrm{ev}_\l(g,t)\}$. 
\end{proof}

\section{Proof of Propositions  \ref{indep} and \ref{zow}}\label{constancy}

We use the notation of \S\ref{zowee}. To define $\tilde f$ explicitly, write $\tilde f(g,t)=((1,g),f_g(t))$, 
where $f_g(t)\colon\Grass(\overline V\oplus V)\to\CCC$ is determined by \S\ref{constr}(ii) and the condition \begin{equation*}f_g(t)(\l\oplus\l):=t(\l)\end{equation*} for any fixed $\l\in\Grass(V)$.

\begin{lemma}\label{embed} The map $\tilde f$ is a homomorphic embedding and is independent of $\l$.
\end{lemma}

\begin{proof}[of Lemma \ref{embed}] According to \S\ref{constr}(ii), $\tilde f$ is independent of $\l$ as long as  
\begin{equation*}
\tau(\l\oplus\l,\l\oplus g\l,\l'\oplus g\l',\l'\oplus\l')=\tau(\l,g\l,g\l',\l'),
\end{equation*}which holds by \S\ref{mi}(iii) and the fact  that $\tau(\l,\l,\l',\l')=0$. Indeed, $\tau(\l,\l,\l',\l')$ is represented by a quadratic space of dimension zero, by Proposition \ref{rankdisc}.

According to \eqref{multiply}, $\tilde f$ is a homomorphism if $c_{(1,g),(1,h)}(\l\oplus\l)=c_{g,h}(\l)$, that is, if
\begin{equation*}\tau(\l\oplus\l,\l\oplus g\l,\l\oplus gh\l)=\tau(\l,g\l,gh\l).
\end{equation*} This again follows from \S\ref{mi}(iii) and Proposition \ref{rankdisc}.

Finally, $\tilde f$ is obviously injective. To show it is continuous, it is enough to notice that for every $\l\in\Grass(V)$ and every $n\geq 0$,  $f$ maps  $N_{\l,n}$ into $N_{\l\oplus\l,n+\dim \l}$; as a result, the subspace topology on $\Mp(V)\subset\Mp(\overline V\oplus V)$ satisfies Proposition \ref{topology}.
\end{proof}

\begin{proof}[of Proposition \ref{zow}]
We have \begin{equation*}
\mathrm{ev}_{\Gamma_1}\circ\tilde f:=f_g(t)(\Gamma_1)=f_g(t)(\l\oplus\l)\cdot\gamma(\tau(\l\oplus\l,\l\oplus g\l,\Gamma_g,\Gamma_1))
\end{equation*}
by \S\ref{constr}(ii). Since $f_g(t)(\l\oplus\l)=t(\l)$, it remains only to check 
\begin{equation}\label{expando}
\tau(\Gamma_g,\Gamma_1,\l\oplus\l)=\tau(\l\oplus\l,\l\oplus g\l,\Gamma_g,\Gamma_1).
\end{equation}
By \S\ref{mi}(i,ii) the difference between the two sides is $\tau(\Gamma_g,\l\oplus\l,\l\oplus g\l)$. This vanishes by Proposition \ref{rankdisc}: it is represented by a quadratic space of dimension zero.
\end{proof}

\begin{proof}[of Proposition \ref{indep}] We note that $f_g(t)(\Gamma_1)$ is independent of $\l$. Moreover, by Proposition \ref{topology}, $(g',t')\mapsto t'(\Gamma_1)$ is locally constant for $g'$ in the open set $N_{\Gamma_1,0}.$ Therefore $(g,t)\mapsto f_g(t)(\Gamma_1)$ is locally constant for $(1,g)\in N_{\Gamma_1,0}$, or, equivalently, for $g\in\Sp(V)''$. 
\end{proof}

\section{Deduction of Theorem 1 from Theorem 2}\label{implication}

  To deduce Theorem 1 from Theorem 2, we have to check
\begin{equation}\label{tbc}\Theta_\l(g,t)\overset?=\pm\gamma(1)^{\dim V-\dim\ker(g-1)-1}\gamma(\det\sigma_g)\end{equation}
with a plus when $F$ is finite or complex and we use the splitting of Proposition \ref{split}. 

According to Proposition \ref{zow}, $\Theta_\l(g,t)=f_g(t)(\Gamma_1)$; but, by \S\ref{constr}(i),  we have $f_g(t)(\Gamma_1)=\pm m_{(1,g)}(\Gamma_1)=\pm m(\Gamma_g,\Gamma_1)$, with plusses when $F$ is finite or complex. Here we must choose the orientations on $\Gamma_1$ and $\Gamma_g$ to be related by the isomorphism $(1,g)$. By definition \eqref{eqm} of $m(\Gamma_g,\Gamma_1)$, it remains to prove

\begin{lemma} With orientations chosen as above, $O_{\Gamma_g,\Gamma_1}=\det\sigma_g\bmod(F^\times)^2$.\end{lemma}

\begin{proof} Consider the isomorphisms
$\Gamma_g \cong V \cong \Gamma_1$ defined by $(x,gx) \mapsto x \mapsto (x,x)$. The hypothesis is just that the orientations of $\Gamma_g$ and $\Gamma_1$ correspond to the same orientation of $V$. These same isomorphisms induce
 \begin{equation*}\Gamma_g/\Gamma_g\cap\Gamma_1\cong V/\ker(g-1)\cong \Gamma_1/\Gamma_g\cap\Gamma_1.\end{equation*}
The symplectic pairing $(\Gamma_g/\Gamma_g\cap\Gamma_1)\otimes(\Gamma_1/\Gamma_g\cap\Gamma_1)\to F$, used to define  $O_{\Gamma_g,\Gamma_1}$ in \S\ref{DA}, induces the pairing $x\otimes y\mapsto\<\sigma_g(y),x\>$ on $V/\ker(g-1)\otimes V/\ker(g-1).$ 
But $\det\sigma_g$ is by its definition in \S\ref{first} the discriminant of this bilinear form.
\end{proof}

\begin{remnum} The proof above is closely related to the following proposition, which is an adaptation of Proposition \ref{rankdisc}.

\begin{prop}\label{rts} There is a quadratic space $(T,q)$ representing $\tau(\Gamma_g,\Gamma_1,\l\oplus \l)$ with rank
\begin{equation*}\dim T=\tfrac12\dim V-\dim\ker(g-1)-\dim g\l\cap \l+2\dim\l\cap\ker(g-1)\end{equation*}
and discriminant
\begin{equation*}\det {q}=(-1)^{\dim\l\cap(g-1)\l}O\gll\det\sigma_g.\end{equation*}
\end{prop}
Formula \eqref{tbc} then follows, alternatively, from Proposition \ref{pd}(iii).
 \end{remnum}

\section{Overview of the Proof of Theorem 2} \label{overview}

The idea used to prove Theorem 2 is to represent the operator $\rho(g,t)$ by an integral kernel and then to compute the trace as the integral along the diagonal (this works quite literally over a finite field).

\subsection{Integral Operators.}
 Recall that for each $\l\in\Grass(V)$ the representation $\rho$ is realized in a Hilbert space completing the space of Schwartz sections $\SSS(\HHH_\l)$ of a certain sheaf (complex line bundle) $\HHH_\l$ on $V/\l$. We recall how to define $\HHH_\l$ in \S\ref{defH}. Let $\HHH_\l^\vee$ be the `dual' sheaf such that a section of $\HHH_\l^\vee\otimes\HHH_\l$ is a measure on $V/\l$.

 We define in \S\ref{secQ} a generalized section $K\lll$ of $\HHH_{\l_1}^\vee\boxtimes\HHH_{\l_2}$ on $V/\l_1\times V/\l_2$. 
In \cite{Th} we proved that the convolution
\begin{equation*}\FF_{21}\colon\SSS(\HHH_{\l_1})\to\SSS(\HHH_{\l_2})\qquad f\mapsto f*K\lll\end{equation*}
 is the operator defined in \cite{LV}, meaning that, for any $\l_1,\ldots,\l_n\in\Grass(V)$, we have
\begin{equation}\label{compos}
\FF_{1n}\circ\FF_{n(n-1)}\circ\cdots\circ\FF_{21}=\gamma(-\tau(\l_1,\ldots,\l_n))\end{equation} and we can realize $\rho$ in the following way.

 Let $\alpha_g:V/\l\times V/\l\to V/g\l\times V/\l,$ $\alpha_g(x,y)=(gx,y). $
Then 
\begin{equation}\label{rho1}
\rho(g,t)\colon\SSS(\HHH_{\l})\to\SSS(\HHH_{\l})\qquad\mbox{is}\qquad f\mapsto t(\l)\cdot  f * \alpha_g^*K\gll.
\end{equation}
\begin{remnum}\label{pullback} The expression $\alpha_g^*K\gll$ should be understood as follows: 
the sheaves $\HHH_{\l_1}^\vee\boxtimes\HHH_{\l_2}$, for varying $\l_1,\l_2$, define a sheaf on \begin{equation*}\{(\l_1,\l_2;w_1,w_2)\,|\,\l_i\in\Grass(V),\,w_i\in V/\l_i\}\end{equation*} with a natural $\Sp(V)\times\Sp(V)$-equivariant structure.  
So $\alpha_g^*K\gll$ is a generalized section of
$\HHH_{\l}^\vee\boxtimes\HHH_{\l}$ on $V/\l\times V/\l$.
\end{remnum}

\subsection{Restriction to the Diagonal.}
 Now let $\Delta\colon V/\l\to V/\l\times V/\l$ be the diagonal, so $\Delta^*(\HHH_{\l}^\vee\boxtimes\HHH_{\l})=\HHH_{\l}^\vee\otimes\HHH_{\l}$ is the sheaf of measures on $V/\l$. Naively, it would follow from \eqref{rho1} that
\begin{equation}\label{bad}\mbox{``}\Tr\rho(g,t)=t(\l)\cdot\int_{V/\l}\Delta^*\alpha_g^*K\gll.\mbox{''}\end{equation} 
This is quite correct when $F$ is finite; let us first consider that case. In doing so, we identify measures and functions on finite sets, using counting measure as a standard.

\subsubsection{}\label{finconc}   We define in \S\ref{theform} (and for any $F$) a quadratic space $(S\gl,q\gl)$ whose class in $W(F)$ is $\tau(\Gamma_g,\Gamma_1,\l\oplus\l)$. By definition, $S\gl$ is a quotient of a certain subspace $\hat S\gl\subset V/\l$ (see \S\ref{hatS}). 

\begin{prop} \label{finrel} Suppose $F$ is finite. Fix $g\in\Sp(V)$. Then
  $\Delta^*\alpha^*K\gll$ is supported on $\hat S\gl$, and, as functions there,
\begin{equation*}
\Delta^*\alpha^*K\gll(x)=\psi(\tfrac12q\gl(x,x))\cdot |F|^{-\tfrac12\dim \l/g\l\cap\l}.
\end{equation*}
\end{prop}
Theorem 2B now follows almost immediately  from \eqref{bad} and Proposition \ref{pg1}---the details, and the proof of Proposition  \ref{finrel}, are worked out in \S\ref{finitefield}.

 \subsubsection{} Now suppose $F$ is infinite. We henceforth restrict ourselves, as we are entitled, to $g$ in the dense open set
\begin{equation}\label{Spl}
\Sp(V)^\l:=\{g\in\Sp(V)''\,|\, g\l\cap \l=0\}
\end{equation}
with $\Sp(V)''$ defined by \eqref{Sp''}. 
From the definition of $K\gll$ (more precisely, from formula \eqref{suppK}) one knows 

\begin{lemma}\label{smooth}
The gene\-ralized sec\-tion $(g,x,y)\mapsto\alpha^*K\gll(x,y)$ is 
smooth when restricted to $\Sp(V)^\l\times V/\l\times V/\l.$
\end{lemma}
 Thus the restriction $\Delta^*\alpha_g^*K\gll$ to the diagonal is a well defined and smooth measure on $V/\l$. 
Let $dq\gl$ be the self-dual measure on $S\gl$. 
   In \S\ref{restrict} we prove the following analogue of Proposition \ref{finrel}.
\begin{prop}\label{rel} Fix $g\in\Sp(V)^\l$. Then $S\gl=V/\l$, and, 
 as measures on $V/\l$,
\begin{equation*}\begin{aligned}
\Delta^*\alpha^*K\gll(x)&=\psi(\tfrac12q\gl(x,x))\cdot\left\|\det(g-1)\right\|^{-1/2}\, dq\gl.\end{aligned}\end{equation*}
\end{prop}

Following \S\ref{finconc}, we should apply the correct infinite versions of \eqref{bad} and Proposition \ref{pg1}. These are given by Lemma \ref{pt} and Proposition \ref{pg}, and 
a little analysis in  \S\ref{proof} completes the proof of Theorem 2A.

\def\OR#1#2{\Omega^\RRR_{#1}\!\(#2\)}
\def\OC#1#2{\Omega_{#1}\!\(#2\)}

\section{Definition of the Sheaf $\HHH_\l$  and the Convolution Kernel $K\lll$}\label{lb}

We must first fix some notation.

\subsection{Conventions on Measures and Densities.}\label{hd}

 For $\alpha\in\RRR$, the space  of $\alpha$-densities on an $F$-vector space $X$ is defined to be the one-dimensional $\RRR$-vector space
\begin{equation*}\OR\alpha{X}=\{\nu\colon\det X\to\RRR\,|\,\nu(\lambda x)=|\lambda|^\alpha \nu(x),\,\forall x\in\det X,\,\lambda\in F\}.\end{equation*}

We identify $\OR{1}{X}$ with the space of real invariant measures on $X$: $\nu\in\OR{1}{X}$ corresponds to the invariant measure that assigns to  $\{a_1v_1+\cdots+a_kv_k\,|\,a_i\in F, |a_i|\leq 1\}$ the volume $\nu(v_1\wedge\ldots\wedge v_k)$, for any basis $v_1,\ldots,v_k$ of $X$. 

\subsubsection{} \label{OR:aut}An isomorphism $f\colon X\to Y$ induces an isomorphism $\OR{\alpha}{X}\to\OR{\alpha}{Y}$, such that $g\in\GL(X)$ acts on $\nu\in\OR{\alpha}{X}$ by $g\cdot\nu=\left\|\det g\right\|^{-\alpha}\nu$, with $\left\|\cdot\right\|$ defined in Theorem 2A.

 We can identify
$\OR{\alpha}{X}\otimes\OR{\beta}{X}=\OR{\alpha+\beta}{X}$
and $\OR{-\alpha}{X}=\OR{\alpha}{X}^*=\OR{\alpha}{X^*}.$
 If $Y\subset X$ then one can identify $\OR{\alpha}{X}=\OR{\alpha}{Y}\otimes\OR\alpha{X/Y}$.

\subsubsection{} Set $\OC\alpha{X}:=\OR\alpha{X}\otimes_\RRR\CCC$. Then \S\ref{OR:aut} works for $\Omega_{\alpha}$ as well as for $\Omega^\RRR_\alpha$.

\subsubsection{}\label{sqrt}
    For $\nu\in\OR{1}{X}$ a positive measure, define $\nu^{1/2}\in\OR{1/2}{X}$ by $\nu^{1/2}(x):=\left|\nu(x)\right|^{1/2}.$ Then $\nu^{1/2}\otimes\nu^{1/2}=\nu$ using $\OR{1/2}{X}\otimes\OR{1/2}{X}=\OR{1}{X}$.

\subsection{Definition of $\HHH_\l$.} \label{defH}
Let $\psi$ be our fixed additive character. For open $U\subset V/\l$ let $\tilde U$ be its pre-image in $V$. Let $\HHH_\l$ be the sheaf on $V/\l$ such that a smooth section $f$ over $U$ is a smooth function $\tilde f\colon \tilde U\to\OC{1/2}{V/\l}$   satisfying the condition
\begin{equation*}\tilde f(v+a)=\psi(\tfrac12\<v,a\>)\cdot \tilde f(v)\qquad\forall v\in V,a\in \l.\end{equation*}

\subsection{Definition of $K\lll$.}\label{gs} \label{secQ} A generalized section $K\lll$ of $\HHH_{\l_1}^\vee\boxtimes\HHH_{\l_2}$
is `the same' as a generalized function $\tilde K\lll\colon V\times V\to
\OC{1/2}{V/\l_1}\otimes \OC{1/2}{V/\l_2}$ satisfying the $\l\times \l$-equivariance condition
\begin{equation*}\tilde K\lll(v+a,w+b)=\psi(\tfrac12\<a,v\>)\cdot \tilde K\lll(v,w)\cdot\psi(\tfrac12\<w,b\>)\end{equation*} for any $a\in \l_1$ and $b\in \l_2$.  

 Our $K\lll$ will be smooth on its support
\begin{equation}\label{suppK}
\supp K\lll=\{(x,y)\in V/\l_1\times V/\l_2\,|\,x-y\in \l_1+\l_2\}.\end{equation}
In other words, 
$\supp \tilde K\lll$ will be the $\l_1\times \l_2$-invariant 
subspace $T\lll\subset V\times V$:
 \begin{equation}\label{eqT}T\lll:=\{(x,y)\in V\times V\,|\,x-y\in \l_1+\l_2\}=\ker\left[V\times V\overset\partial\too V/(\l_1+\l_2)\right].
\end{equation}
$\tilde K\lll$ is constructed from the following quadratic form $Q\lll$ on $T\lll$.

\subsubsection{}\label{defQ}

Given $(x,y)\in T\lll$ write $x-y=a_1+a_2$ with $a_i\in \l_i$. Then
\begin{equation}\label{Q}Q\lll((x,y)):=\<a_1,x\>+\<a_2,y\>.\end{equation}

\begin{lemma}\label{lemQ} 
\begin{enumerate}
\item[(i)] $Q\lll((x,y))$ is independent of the choice of $a_1,a_2$.
\item[(ii)] For $b_i\in \l_i$, $Q\lll((x+b_1,y+b_2))=Q\lll((x,y))+\<b_1,x\>+\<y,b_2\>.$
\end{enumerate}
\end{lemma}

\begin{proof} For (i), if $x-y=(a_1+\epsilon)+(a_2-\epsilon)$, with $\epsilon\in \l_1\cap \l_2$, then
\begin{equation*}\<a_1+\epsilon,x\>+\<a_2-\epsilon,y\>=\<a_1,x\>+\<a_2,y\>+\<\epsilon,x-y\>\end{equation*}
but $\<\epsilon,x-y\>=0$ since $x-y\in\l_1+\l_2$.  The second statement is similarly easy.
\end{proof}

\subsubsection{Definition of $\tilde K\lll$.} \label{defK}

 Let $\delta\lll$ be the extension-by-zero of the constant function $1$ from $T\lll$ to $V\times V$. It is naturally a generalized function on $V\times V$ with values in $\OC{1}{(V\times V)/T\lll}=\OC{1}{V/(\l_1+\l_2)}$.

Let $\mu\lll$ be the element of 
$\OC{1}{(\l_1+\l_2)/\l_1\cap\l_2}$
corresponding to self-dual measure on the symplectic space $(\l_1+\l_2)/\l_1\cap \l_2$. Let $\mu\lll^{1/2}$ be its square root (\S\ref{sqrt}).

Set 
\begin{equation*}\tilde K\lll(v,w):=\psi(\tfrac12Q\lll((v,w)))\cdot \delta\lll(v,w)\cdot\mu\lll^{1/2}.\end{equation*}
Since 
$\OC{1}{(\l_1+\l_2)/\l_1\cap\l_2}=\OC{1}{(\l_1+\l_2)/\l_2}\otimes\OC{1}{(\l_1+\l_2)/\l_1},$ one sees that $\tilde K\lll$ is a generalized function on $V\times V$ with values in $\OC{1/2}{V/\l_1}\otimes\OC{1/2}{V/\l_2}$.  
By Lemma \ref{lemQ}(ii) and the first sentence of \S\ref{gs}, $\tilde K\lll$ determines a generalized section $K\lll$ of $\HHH_{\l_1}^\vee\boxtimes\HHH_{\l_2}$. 

\subsubsection{Definition of $\eta\lll$.}\label{eta} For technical purposes in \S\ref{restrict} we will also need:
\begin{equation*}\eta\lll(v,w):=\delta\lll(v,w)\cdot\mu\lll^{1/2}.\end{equation*}
Since $\eta\lll$ depends on $(v,w)$ only modulo $\l_1\times \l_2$, we may consider it as a generalized function on $V/\l_1\times V/\l_2$, supported on \eqref{suppK}, with values in $\OC{1/2}{V/\l_1}\otimes\OC{1/2}{V/\l_2}$. The pull-back $\alpha_g^*\eta\gll$ should be understood much as in Remark \ref{pullback}: it is a generalized function on $V/\l\times V/\l$ with values in $\OC{1}{V/\l}$.

\section{Two Quadratic Spaces}\label{theform}

Now we define the quadratic spaces $(S\gl,q\gl)$ and $(S'\gl,q\gl)$ mentioned in \S\ref{defq1} and \S\ref{finconc} (they are dual, hence isometrically isomorphic to one another). We first give an abstract definition, and then explicit formulas in \S\ref{expform}. In \S\ref{MaslovInterp}, we show that both forms represent $\tau(\Gamma_g,\Gamma_1,\l\oplus\l)$ in $W(F)$. Finally, in \S\ref{maktouf} we detail the connection to Maktouf's work \cite{Ma} mentioned in \S\ref{mak1}.

\subsection{Definition of $S\gl$ and $S'\gl$.}\label{Dq} Let $A,B,C,D$ denote the rows (they are complexes) in the commutative diagram
\begin{equation*}
\xymatrix@R=0.3in@C=.3in{
& g\l\cap\l \ar[r]^{1-g\inv}\ar[d]^{g\inv}\ar[d]
&  {\l} \ar[r]\ar[d]
&  V/(g-1)V\ar@{=}[d]
\\
&  {\l} \ar[r]^{g-1}\ar[d]
&  g\l+\l \ar[r]\ar[d]
&  V/(g-1)V\ar@{=}[d]
\\
{\ker(g-1)}\ar[r]\ar@{=}[d]
&  V \ar[r]^{g-1}\ar[d]
&  V \ar[r]\ar[d]
&  V/(g-1)V
\\
{\ker(g-1)}\ar[r]
&  V/\l \ar[r]^{g-1}
&  V/(g\l+\l)
}
\end{equation*}

\begin{definition} Let $S\gl$ be the cohomology of $D$ at its center term, and $S'\gl$ be the cohomology of $A$ at its center term. Then $S\gl$ and $S'\gl$ are dual to one another under the symplectic pairing, because indeed $A$ and $D$ themselves are dual. 
\end{definition}

\subsection{Definition of $q\gl$ and $q'\gl$.} \label{DPhi}

Let $W$ denote the cohomology of the complex $B$ at its center term. The map $A\to B$ is a quasi-isomorphism of complexes, giving an isomorphism $S'\gl\to W$. $B\to C\to D$ is a short exact sequence of complexes, and $C$ is acyclic, so the boundary map gives an isomorphism $S\gl\to W$ on cohomology. Together we have an isomorphism
\begin{equation*}\Phi\gl\colon S\gl\to W\leftarrow S'\gl.\end{equation*}

\begin{definition} Set 
\begin{equation*}q\gl(x,y):=\<\Phi\gl(x),y\>\qquad q'\gl(a,b):=\<a,\Phi\gl\inv(b)\>\end{equation*}
 for all $x,y\in S\gl$, $a,b\in S'\gl.$ In particular, $\Phi\gl\colon S\gl\to S'\gl$ is an isometric isomorphism. 
\end{definition}

\subsection{Explicit Formulas.}\label{expform} Let us now give explicit formulas for $q\gl$ and $q'\gl$ and show that they are symmetric forms.

\subsubsection{}\label{hatS} First observe that $S\gl$ is a quotient of
\begin{equation}\label{eqS}\hat S\gl:=\{x\in V/\l\,|\,(g-1)x\in g\l+\l\}\subset V/\l\end{equation}
and $S'\gl$ is a quotient of
\begin{equation}\label{eqS'}\hat S'\gl:=\l\cap(g-1)V\subset \l.\end{equation}
We will really give formulas for $q\gl$ and $q'\gl$ pulled back to $\hat S\gl$ and $\hat S'\gl$.

\begin{remnum}\label{remS} For $g\in\Sp(V)''$, as in \eqref{Sp''}, $S\gl=\hat S\gl$ and $\hat S'\gl=\l$. Moreover, if $g\in\Sp(V)^\l$, as in \eqref{Spl}, then $S\gl=\hat S\gl=V/\l$ and $S'\gl=\hat S'\gl=\l$.
\end{remnum}

\subsubsection{}Given $x,y\in \hat S\gl$, there exist $a,b\in\l$ with $(g-1)x\equiv(ga+b)\bmod (g-1)\l$.  Then
\begin{equation} \label{expq}q\gl(x,y)=\<a+b,y\>.\end{equation}
(Indeed, $\Phi\gl(x)=a+b$ because $a+b\equiv ga+b\bmod (g-1)\l$.)

  Given $a,b\in \hat S'\gl$, suppose $b=(g-1)y$. Then
\begin{equation}\label{expq'}q'\gl(a,b)=\<a,y\>.\end{equation}

\begin{prop}\label{aresym} Both $q\gl$ and $q'\gl$ are symmetric forms.
\end{prop}
\begin{proof} 
The two forms are isometric, so it suffices to consider $q'\gl$.
 Given $a,b\in \hat S'\gl$, suppose $a=(g-1)x, b=(g-1)y$. Then
\begin{equation*}q'\gl(a,b)=\<(g-1)x,y\>=\<x,(g\inv-1)y\>=\<x,-g\inv b\>=\<gx,-b\>.\end{equation*}
Now  $gx-x=a\in\l$, so $\<gx,-b\>=\<x,-b\>=\<b,x\>=q'\gl(b,a).$
\end{proof}

\subsection{Interpretation via the Maslov Index.}\label{MaslovInterp}

\begin{prop}\label{MI} The class of $(S\gl,q\gl)$ and therefore of $(S'\gl,q'\gl)$ in $W(F)$ equals $\tau(\Gamma_g,\Gamma_1,\l\oplus\l)$. In fact, these quadratic spaces satisfy Proposition \ref{rts}.
\end{prop}
\begin{proof}
Recall from \cite{Th} that $\tau(\Gamma_g,\Gamma_1,\l\oplus \l)$ is represented by a degenerate quadratic form $q$ on
\begin{equation*}\hat T:=\{(x,y,z)\in\Gamma_g\times\Gamma_1\times(\l\oplus \l)\,|\,x+y+z=0\in\overline V\oplus V\}.\end{equation*}
Namely, $q((x,y,z))=\<x,z\>$, where we pair
using the symplectic form on $\overline V\oplus V$.

We claim that the map $f\colon \hat T\to \hat S'\gl$, \begin{equation*}f\colon((x,gx),(y,y),(a,b))\mapsto a-b\end{equation*}
 descends to an isometric isomorphism on the nondegenerate quotients.

 First we observe that the condition $s:=((x,gx),(y,y),(a,b))\in \hat T$ implies that $a-b=(g-1)x$, so $f$ does have the right target. It is also an isometry since
\begin{equation*}q(s)=\<gx,b\>-\<x,a\>=\<x,b\>-\<x,a\>=\<a-b,x\>=q'\gl(f(s),f(s)).\end{equation*}
The second equality holds because $gx\equiv x\bmod\l$, the fourth from \eqref{expq'}.

Any isometry descends to an injective map of nondegenerate quotients, so
it remains to observe that the nondegenerate quotient $S'\gl$ of $\hat S'\gl$, defined in \S\ref{Dq},  has the same dimension as the nondegenerate quotient $T$ of $\hat T$ (cf. Proposition \ref{rts}), namely,
\begin{equation}\label{dim}\begin{aligned}
\dim S'\gl=\dim S\gl=\tfrac12&\dim V-\dim\ker(g-1) \\
& -\dim g\l\cap\l+2\dim{\l\cap \ker(g-1)}.\end{aligned}\end{equation}
\end{proof}

\subsection{Relation with Maktouf's Construction.}\label{maktouf}

 Let $(s,t_s)$ be the semisimple part of $(g,t)\in\Mp(V)$. Assume for simplicity that $(s-1)$ is invertible. Maktouf \cite{Ma} asserts that there exists an $s$-stable symplectic decomposition $V=V_1\oplus V_2$ and a Lagrangian $\l_i\in V_i$ such that $s\l_1=\l_1$ and $s\l_2\cap \l_2=0$. He then defines a quadratic form $q_{\mathrm{Mak}}$ on $\l_2$ by $q(a,a)=\<a,(s\inv-1)\inv a\>$ and sets
\begin{equation*}\Phi(g,t):=t_s(\l_1\oplus\l_2)\cdot\gamma(-q_{\mathrm{Mak}}).\end{equation*}

Suppose that $g=s$ is semisimple and $g-1$ is invertible. 
Choose $\l:=\l_1\oplus\l_2\subset V$, and consider $q'\gl$ as a form on $\hat S'\gl=\l$.  We will show that $q_{\mathrm{Mak}}=-q'\gl$, so $\Phi=\Theta_\l$.

\begin{lemma}\label{kerq2} We have $\l_1=\ker q'\gl$. Thus $q'\gl$ defines a nondegenerate bilinear form on $\l_2$.
\end{lemma}

\begin{proof} Clearly $g\l\cap \l=\l_1$, and $\l_1$ is $(1-g\inv)$-stable. But $q'\gl$ was defined to be nondegenerate on $S'\gl=\l/((1-g\inv)(g\l\cap\l))$.
\end{proof}

\begin{prop}\label{qq} We have $q'\gl=-q_{\mathrm{Mak}}$ on $\l_2$. \end{prop}
\begin{proof*}
Suppose $a\in \l_2$. Then 
\multbox\begin{eqnarray*}
q_{\mathrm{Mak}}(a,a)&&=\<a,(g\inv-1)\inv a\>=\<(g-1)(g-1)\inv a,(g\inv-1)\inv a\>\\
&&=\<(g-1)\inv a,(g\inv-1)(g\inv-1)\inv a\>\\&&=\<(g-1)\inv a,a\>=-q'\gl(a,a).
\end{eqnarray*}
\end{proof*}

\section{Proof of Proposition \ref{finrel} and Theorem 2B}\label{finitefield}

Let $F$ be a finite field; we again identify measures and functions on finite sets.

\subsection{Proof of Proposition \ref{finrel}.}
\label{finrestrict}
For any $g\in\Sp(V)$, consider the composition 
\begin{equation*}\label{dag}
\begin{CD}
V/\l@>\Delta>> V/\l\times V/\l @>{\alpha_g}>> V/g\l\times V/\l\end{CD}\qquad x\mapsto(x,x)\mapsto(gx, x).
\end{equation*}
Proposition \ref{finrel} amounts to the following lemma.

\begin{lemma}\label{finlrestrict} Fix any $g\in\Sp(V)$.
\begin{enumerate}
\item[(i)] $\Delta\inv\circ\alpha_g\inv(\supp K\gll)=\hat S\gl\subset V/\l.$
\item[(ii)] $\Delta^*\alpha_g^*Q\gll(x)=q\gl(x,x)$ for any $x\in \hat S\gl$. 
\item[(iii)] $\Delta^*\alpha^*\eta\gll(x)=|F|^{-\tfrac12\dim{\l/g\l\cap\l}}$ for any $x\in \hat S\gl$. 
\end{enumerate}
(Notation: see \S\ref{defQ} for $Q\gll$ and  \S\ref{eta} for $\eta\gll$.)
\end{lemma}

\begin{proof}
The first statement follows from \eqref{suppK}, \eqref{eqS}.  
The second statement follows from \eqref{expq} and \S\ref{defQ}.
For the third, recall that $\eta\gll$ is defined in \S\ref{eta} to be constant on its support (which is $\hat S\gl$ by part (i)), with value 
\begin{equation*}\mu\gll^{1/2}=\(|F|^{-\tfrac12\dim\((g\l+\l)/g\l\cap\l\)}\)^{1/2}= |F|^{-\tfrac12\dim \l/g\l\cap\l}.\end{equation*}
Here we use that if $A$ has a nondegenerate bilinear form, then the self-dual measure on $A$ is $|F|^{-\dim A/2}$.
\end{proof}

\subsection{Proof of Theorem 2B.}
Using \eqref{rho1} and Proposition \ref{finrel} above, we have
\begin{equation*}
\Tr\rho(g,t)=t(\l)\cdot\sum_{x\in \hat S\gl} \psi(\tfrac12q\gl(x,x))\cdot |F|^{-\tfrac12\dim \l/g\l\cap\l}.
\end{equation*}
Applying Propositions \ref{pg1} and \ref{MI}, we obtain
\begin{equation}\begin{aligned}
\label{almost}\Tr\rho(g,t)&=t(\l)\cdot |F|^{\tfrac12\(\dim S\gl-\dim \l/g\l\cap\l\)+\dim\ker q\gl}\cdot \gamma(q\gl)\\
&=|F|^{\tfrac12\(\dim S\gl-\dim \l/g\l\cap\l\)+\dim\ker q\gl}\cdot\Theta_\l(g,t)
\end{aligned}
\end{equation}
where $\ker q\gl$ is the kernel of $q\gl$ as a form on $\hat S\gl$.
Now, by construction of $S\gl$ in  \S\ref{Dq}, we have
\begin{equation*}
\dim \ker q\gl=\dim\ker(g-1)-\dim\l\cap\ker(g-1)
\end{equation*}
as well as formula \eqref{dim} for $\dim S\gl$. 
This with \eqref{almost} establishes Theorem 2B.

\section{Proof of Proposition \ref{rel} and Theorem 2A} 

\subsection{Proof of Proposition \ref{rel}.} \label{diagonal}

\label{restrict}
Proposition \ref{rel} amounts to Lemma \ref{finlrestrict}(i,ii), which still hold, Remark \ref{remS}, and the following infinite version of Lemma \ref{finlrestrict}(iii).

\begin{lemma}\label{lrestrict} 
 If $g\in\Sp(V)^\l$ then $\Delta^*\alpha_g^*\eta\gll=\left\|\det (g-1)\right\|^{-1/2}dq\gl$ as measures on $V/\l$. 
\end{lemma}

\begin{proof}
 It is clear that $\nu_\l:=\Delta^*\alpha_g^*\eta\gll$ is an invariant measure on $\hat S\gl=S\gl=V/\l$; that is, $\nu_\l\in\OC{1}{V/\l}$ in the notation of \S\ref{hd}. 
Let $\Phi\gl\colon V/\l \to \l$  be the isomorphism defining $q\gl$, as in \S\ref{DPhi}.   The claim is that 
\begin{equation}\label{muit1}\<\nu_\l,{\Phi\gl}_*\nu_\l\>=\left\|\det (g-1)\right\|\inv \end{equation} 
under the pairing of $\Omega_1(V/\l)$ with $\Omega_1(\l)$ induced by the symplectic form. 

Let $\mu_V$ be the self-dual measure on $V$. Fix $\omega\in\det \l$ and consider $g\omega\in\det(V/\l)$. It is easy to see that $\nu_\l$, as a function $\det(V/\l)\to\RRR$ in the sense of \S\ref{hd}, satisfies 
\begin{equation*}\nu_\l(g\omega)=\left(\mu_V(g\omega\wedge\omega)\right)^{1/2}.\end{equation*}
  On the other hand $\Phi\gl\inv(x)=(g-1)\inv x\bmod \l$ for any $x\in\l$; thus 
\begin{equation*}
{\Phi\gl}_*\nu_\l(\omega)=\nu_\l((g-1)\inv\omega)=
\nu_\l(g\omega)\cdot\frac{\mu_V((g-1)\inv\omega\wedge\omega)}{\mu_V(g\omega\wedge\omega)},
\end{equation*}
the fraction here being the absolute value of the ratio of $(g-1)\inv\omega$ and $g\omega$ as elements of $\det(V/\l)$. 
 Since $(g-1)\inv\omega\wedge\omega=(g-1)\inv(\omega\wedge(g-1)\omega)=(g-1)\inv(\omega\wedge g\omega)$, we have all together
\begin{equation*}
{\Phi\gl}_*\nu_\l(\omega)=\left\|\det (g-1)\right\|\inv\left(\mu_V(g\omega\wedge \omega)\right)^{1/2}.\end{equation*}
Since  $\<\nu_\l,{\Phi\gl}_*\nu_\l\>=\frac{\nu_\l(g\omega)\cdot{\Phi\gl}_*\nu_\l(\omega)}{\mu_V(g\omega\wedge\omega)}$,
we obtain \eqref{muit1}.
\end{proof}

\subsection{Proof of Theorem 2A.}\label{proof}

\subsubsection{}\label{char} First recall how $\Tr\rho$ is defined as a generalized function on $\Mp(V)$. If $M$ is a compactly supported smooth function on $\Mp(V)$ then 
\begin{equation*}\rho(M)\colon f\mapsto\int_{(g,t)\in\Mp(V)} M(g,t)\rho(g,t)f\end{equation*}
is a trace-class operator on $\SSS(\HHH_{\l_1})$, and
\begin{equation*}\<\Tr\rho,M\>:=\Tr\rho(M).\end{equation*} 

\subsubsection{} \label{Mpl} Assume from now on that $M$ is supported over $\Sp(V)^\l$, see \eqref{Spl}. Here is a valid version of \eqref{bad}.

\begin{lemma}\label{pt}
$\displaystyle{{\Tr\rho(M)=\int_{x\in V/\l}\int_{(g,t)\in\Mp(V)} \!\!\!\!M(g,t)\cdot t(\l)\cdot \Delta^*\alpha_g^*K\gll(x).}}$
\end{lemma}

\begin{proof} As in \eqref{rho1}, the operator $\rho(M)$ is represented by the integral kernel
\begin{equation}\label{rhoM}(x,y)\mapsto\int_{(g,t)\in\Mp(V)} M(g,t)\cdot t(\l)\cdot \alpha_g^*K\gll(x,y)\end{equation}
on $V/\l\times V/\l$.  We know $\rho(M)$ is trace class. For its trace to be the integral of \eqref{rhoM} along the diagonal, it suffices that \eqref{rhoM} be smooth.   This is clear from Lemma \ref{smooth}.  
\end{proof}

\subsubsection{}

 Now choose a function $h$ on $A:=V/\l$ as in \S\ref{Dirac}, and again set $h_s(x):=h(sx)$. We claim that
\begin{equation}\label{lim}\Tr\rho(M)=\lim_{s\to0}\Tr(h_s\cdot\rho_M).\end{equation}
Indeed, adapting Lemma \ref{pt}, we have 
\begin{equation}\label{e}\Tr (h_s\cdot\rho(M))=\int_{x\in V/\l}\int_{(g,t)\in \Mp(V)} M(g,t)\cdot h_s(x)\cdot t(\l)\cdot \Delta^*\alpha_g^*K\gll(x).\end{equation}
The outer integral converges absolutely and uniformly in $s$,  being dominated by $\Tr\rho(M)$.  Therefore
\begin{equation*}\lim_{s\to 0} \Tr (h_s\cdot\rho(M))=\int_{x\in V/\l}\int_{(g,t)\in\Mp(V)} M(g,t)\cdot t(\l)\cdot \Delta^*\alpha_g^*K\gll(x)=\Tr \rho(M)\end{equation*}
as desired.

 \subsubsection{} On the other hand, we may exchange the order of integration in \eqref{e}. 
\begin{equation}\label{f}\Tr (h_s\cdot\rho(M))= \int_{(g,t)\in\Mp(V)} M(g,t)\cdot\int_{x\in V/\l} h_s(x)\cdot t(\l)\cdot \Delta^*\alpha_g^*K\gll(x).
\end{equation}
Now, from Propositions \ref{rel} and \ref{pg} and Remark \ref{remS}, we find, for $g\in\Sp(V)^\l$,
\begin{equation*}\lim_{s\to0} \int_{x\in V/\l} h_s(x)\cdot\Delta^*\alpha_g^*K\gll(x) =\left\|\det(g-1)\right\|^{-1/2}\cdot \gamma(q\gl).\end{equation*}
Moreover, due to Lemma \ref{smooth}, this limit converges uniformly for $g$ in a compact set.  Uniformity allows us to calculate $\lim_{s\to 0}\Tr(h_s\cdot\rho(M))$ from \eqref{f} by taking the limit $s\to 0$ inside the first integral to obtain 
\begin{equation*}\lim_{s\to0}\Tr(h_s\cdot\rho(M))=\int_{(g,t)\in\Mp(V)} M(g,t)\cdot t(\l)\cdot\left\|\det(g-1)\right\|^{-1/2}\cdot \gamma(q\gl).\end{equation*}
This combined with \eqref{lim} and Proposition \ref{MI} completes the proof of Theorem 2A.

\affiliationone{%
   Teruji Thomas\\
   Merton College\\
   Oxford OX1 4JD\\
   United Kingdom\\
   \email{jtthomas@uchicago.edu}}

\begin{thebibliography}{19}
\bibitem{Ad} {\bibname J. Adams}, `Character of the oscillator representation', {\em Israel J. Math.} {98} (1997) 229--252.
\bibitem{Bu} {\bibname D. Bump}, {\em Automorphic forms and representations}, Cambridge Studies in Advanced Mathematics 55 (Cambridge University Press, 1998).
\bibitem{GH} {\bibname S. Gurevich and R. Hadani}, `The geometric Weil representation', to appear. Preprint at www.arxiv.org/math.RT/0610818.
\bibitem{Howe} {\bibname R. Howe}, `On the character of Weil's representation', {\em Trans. Amer. Math. Soc.} {177} (1973) 287--298.
\bibitem{KS} {\bibname M. Kashiwara and P. Schapira}, {\em Sheaves on Manifolds}, Grundlehren der mathematischen Wissenschaften 292 (Springer, Berlin, 1990).
\bibitem{Lam} {\bibname T.Y. Lam}, {\em The algebraic theory of quadratic forms} (W.A. Benjamin, Reading, 
Mass., 1973).
\bibitem{LV} {\bibname G. Lion and M. Vergne}, {\em The Weil representation, Maslov index and Theta series}, Progress in Mathematics 6 (Birkha\" user, Boston, 1980).
\bibitem{LP} {\bibname G. Lion and P. Perrin}, `Extension des representations de groupes unipotents $p$-adiques: calculs d'obstructions', {\em Non Commutative Harmonic Analysis and Lie Groups (Marseille-Luminy, 1980)} (eds J. Carmona and M. Vergne), Lecture Notes in Mathematics 880 (Springer, Berlin, 1981), pp. 337--356.  
\bibitem{Ma} {\bibname K. Maktouf}, `Le caract\`ere de la repr\'esentation m\'etaplectique et la formule du caract\`ere pour certaines repr\'esentations d'un groupe de Lie presque alg\'ebrique sur un corps $p$-adique', {\em J. Functional Analysis} {164} (1999) 249--339.
\bibitem{Pe} {\bibname P. Perrin}, `Repr\'esentations de Schr\"odinger, indice de Maslov et groupe metaplectique', {\em Non Commutative Harmonic Analysis and Lie Groups (Marseille-Luminy, 1980)} (eds J. Carmona and M. Vergne), Lecture Notes in Mathematics 880 (Springer, Berlin, 1981), pp. 370--407.
\bibitem{Th} {\bibname T. Thomas}, `The Maslov index as a quadratic space.' {\em Math. Res. Lett.} {13} (2006) 985--999. Expanded electronic version at www.arxiv.org/math.SG/0505561/.
\bibitem{Torasso} {\bibname P. Torasso}, `Sur le caract\`ere de la repr\'esentation de Shale-Weil de $\Mp(n,\RRR)$ et $\Sp(n,\CCC)$', {\em Math. Ann.} {252} (1980) 53--86.
\bibitem{We} {\bibname A. Weil}, `Sur certains groupes d'op\'erateurs unitaires', {\em Acta Math.} {111} (1964) 143--211.
\end{thebibliography}
\end{document}